\documentstyle[amscd,righttag,amssymb]{amsart}
\newmathalphabet*{\script}{eus}{m}{n}
\addtoversion{bold}{\script}{eus}{b}{n}
\newtheorem{theorem}{Theorem}[section]

\newtheorem{lemma}[theorem]{Lemma}
\newtheorem{corollary}[theorem]{Corollary}

\theoremstyle{definition}

\theoremstyle{remark}
\newtheorem{remark}[theorem]{Remark}
\numberwithin{equation}{section}

\newcommand{\ad}{\operatorname{Ad}}
\newcommand{\AAA}{{\Bbb A}}
\newcommand{\CC}{{\Bbb C}}

\newcommand{\PP}{{\Bbb P}}

\newcommand{\SL}{\operatorname{SL}}
\newcommand{\sln}{\SL_{n}(\CC)}
\newcommand{\liesln}{\frak{s}\frak{l}_n(\CC)}
\newcommand{\cmp}[3]{Comm. Math. Phys. {\bf #1} (#2), #3}
\newcommand{\pl}[3]{Phys. Lett. {\bf B #1} (#2) #3}
\newcommand{\np}[3]{Nucl. Phys. {\bf B #1} (#2), #3}

\newcommand{\ijmp}[3]{Int. Jour. Mod. Phys. A{\bf #1} (#2), #3}

\newcommand{\tams}[3]{Trans. AMS {\bf #1} (#2), #3}
\newcommand{\jgp}[3]{Jour. Geom. and Phys. {\bf #1} (#2), #3}
\newcommand{\dmj}[3]{Duke Math. J. {\bf #1} (#2), #3}
\newcommand{\hepth}[1]{{\tt hep-th}/#1}
\newcommand{\qalg}[1]{{\tt q-alg}/#1}

\begin{document}
\begin{flushright}Matematisk Institut, Aarhus Universitet\\
Preprint series 1995 No. 15
\\[7pt]
  {\tt q-alg/9508003}
\end{flushright}
\title[Miura transformation]{Toda Fields on Riemann
  Surfaces: remarks on the Miura Transformation}
\author{E. Aldrovandi}
\address{Department of Mathematics\\SUNY at Stony
Brook\\Stony Brook, NY 11794-3651\\ USA}
\email{ettore@@math.sunysb.edu}
\thanks{Supported by Consiglio Nazionale delle Ricerche,
CNR-NATO grant}
\maketitle
\begin{abstract}
  We point out that the Miura transformation is related to a
  holomorphic foliation in a relative flag manifold over a
  Riemann Surface. Certain differential operators
  corresponding to a free field description of $W$--algebras
  are thus interpreted as partial connections associated to
  the foliation.
\end{abstract}

\section{Introduction and background}
The Miura transformation plays an important role in the
free--field approach to the theory of integrable systems.
This is very effective in the case of a cylindrical
topology, where it has been shown how the correspondence
between solutions to the Toda field equations and
free--fields is basically one to one~\cite{ABBP}, thereby
concluding a program initiated by Gervais and
Neveu~\cite{GeNe}. Free fields are also a powerful tool in
describing the relevant $W$--algebra structures and in
addressing quantization problems~\cite{FaLy,BiGe,GeNe}.

However, the above mentioned correspondence becomes
problematic as soon as one tries to directly extend it to
Riemann Surfaces of higher genus. A ``no--go theorem'' has
been devised in ref.~\cite{AlBo} for the uniformization
solution to the Liouville equation, and the same analysis
applies to the geometric approach to the Toda equations
pursued in ref.~\cite{AlFa}. These results prompt for an
investigation of the meaning of the Miura transformation in
presence of a non trivial topology.

Our aim will be to show in what sense
the well known local calculations
to be recalled below can be recast in a global setting and
to exhibit an interpretation for the relations defining the
Miura transformation. This should be of interest in
connection with the problem of determining to what extent
free fields do describe Toda Field theories on a Riemann
surface.

\subsection{}\label{local-1}
It is well known that classical $W$-algebras can be
described in terms of gauge equivalence classes of chiral
connections under lower triangular gauge
transformations~\cite{ItZu}. The customary $W$-fields appear
when one fixes the gauge to
\begin{equation}
\nabla = \partial +
\begin{pmatrix}
0&1&0&\dots &0\\
0&0&1&\dots &0\\
\hdotsfor{4}&1\\
w_n&w_{n-1}&\dots&\dots&0
\end{pmatrix} dz ,\quad T\equiv w_2
\label{W-form}
\end{equation}
hereafter called ``Drinfel'd-Sokolov'' gauge. (We
follow~\cite{AlFa,BoX,ZuDS}, a different naming convention
was adopted in~\cite{AlBo}.) This is true in particular for
the zero curvature representation of the Toda equations in
the holomorphic gauge~\cite{AlFa,ZuDS}. On the other hand,
there are other relevant gauges~\cite{BoX}, and in
particular the ``diagonal gauge''
\begin{equation}
\widetilde\nabla = \partial +
\begin{pmatrix}
p_1&1&0&\dots&0\\
0&p_2&1&\dots&0\\
\hdotsfor{4}&1\\
0&\dots&\dots&\dots&p_n
\end{pmatrix} dz ,\quad p_1 +\dots +p_n=0
\label{DS-form}
\end{equation}
gives the above mentioned free fields.
It is customary to introduce new
variables $a_1,\dots ,a_{n-1}$ through
\begin{align}
p_1 &= a_1\notag\\
p_k &= a_k - a_{k-1}\, ,\quad k=2,\dots ,n-1\notag\\
p_n &= -a_{n-1}\notag
\end{align}
in order to enforce the null trace condition in
\eqref{DS-form}. One can connect the two slices
\eqref{W-form} and \eqref{DS-form} by a lower triangular
gauge transformation
\begin{equation}
\nabla \mapsto \widetilde\nabla =
\bold{N}_-^{-1}\circ\nabla\circ \bold{N}_-
\label{gauge}
\end{equation}
with
\begin{equation}
\bold{N}_- =\begin{pmatrix}
            1&0&0&\dots&0\\
            a_1&1&0&\dots&0\\
            \ast &a_2&1&\dots&0\\
            \hdotsfor{5}\\
            \ast &\dots&\ast &a_{n-1}&1
            \end{pmatrix}
\end{equation}
(the other terms are determined in terms of $a_1,\dots
,a_{n-1}$) which yields a relation
\begin{equation}
w_i=w_i[\{a_j\},\{\partial a_j\},\dots
,\{{\partial}^{i-1}a_j\}]
\label{miura}
\end{equation}
expressing the $W$-fields in terms of the free fields. This
is precisely the {\em Miura transformation.\/} For example,
in the rank two case one gets
\begin{equation}\label{miura-2}
T(z) = - \partial_z a(z) + a(z)^2\, ,
\end{equation}
while if $n=3$
\begin{align}
w_2 &= -\partial_z a_1 +a_1^2 -\partial_z a_2 + a_2^2 - a_1
a_2\nonumber\\
w_3 &= \partial_z (\partial_z a_1 - a_1^2) +
a_1 (\partial_z a_2 - a_2^2) + a_1^2 a_2\, ,\label{miura-3}
\end{align}
and so on.

\subsection{}\label{local-2}
One is typically interested in a local basis of flat
sections
\begin{equation}\label{fundam-W}
\nabla \bold{M} =0\, .
\end{equation}
The preceding gauge transformation is encoded in the (Gauss)
factorization
\begin{equation}
\bold{M}=\bold{N}_-\, \bold{B}_+
\label{gauss}
\end{equation}
where $\bold{N}_-$ is a lower unipotent matrix and
$\bold{B}_+$ is an upper triangular matrix. This splitting
yields a corresponding basis $\bold{B}_+$ of flat sections
for $\widetilde\nabla$ and automatically provides a solution
to \eqref{miura}. This is the way the Miura transfomation is
hidden in the free field approach to Toda Field theory
exploited in~\cite{ABBP}.
\begin{remark}
The factorization \eqref{gauss} implies taking ratios among
various entries in the matrix $\bold{M}$, thus it is bound
to be ill--defined at certain points. This problem is
usually dealt with by assuming the Gauss factorization can
be performed. We shall characterize the locus where this can
actually be done.
\end{remark}

\subsection{}
The two ${\cal D}$-modules \eqref{W-form} and
\eqref{DS-form} can actually be defined in a more general
setting. (See, in this connection, ref. \cite{Fre95}.)
 Let $C$ be a compact Riemann surface of genus $g>1$. It has
been shown that \eqref{W-form} defines a holomorphic
connection on the vector bundle
${E}=J^{n-1}({{K}_C}^{-\frac{n-1}{2}})$ of
$(n-1)$-jets of holomorphic sections of
${{K}_C}^{-\frac{n-1}{2}}$, where ${K}_C$ is the
canonical line bundle~\cite{AlFa,Zucc}. On the other hand,
\eqref{DS-form} defines a {\em meromorphic\/} connection on
the vector bundle
$\widetilde{E}=\bigoplus_{r=0}^{n-1}{{K}_C}^{-\frac{n-1}{2}+r}$
\cite{AlBo,ZuDS}.
(Because of Weil's theorem, $\widetilde{E}$ does not admit
analytic connections~\cite{Atiyah,Gun67}.) By the monodromy
analysis in~\cite{AlBo}, which easily carries over to the
rank $n$ case, there are no (meromorphic) morphisms from the
pair $(E,\nabla )$ to $(\widetilde{E},\widetilde\nabla ).$
This is especially true when the former results from the
zero--curvature representation of the Toda
equations. Therefore the local procedure leading to the
gauge transformation cannot be globally implemented, so that
on a Riemann Surface we have two radically distinct objects.

The above observation questions the idea of extending the
approach to classical $W$-algebras via gauge equivalence
classes of connections to a higher genus surface. Thus the
line we shall stick to in the following will be to consider
the vector bundle $E$ and the analytic connection
\eqref{W-form} as the definition for $W$-type objects on a
Riemann Surface. This is the same one as in ref.~\cite{AlFa}
and it is in line with other geometric approaches, such
as~\cite{suresh}. (Note that this is also the opposite
to~\cite{ZuDS}.)

\subsection{}
Let us describe our result in some more detail. (See however
the next section for the exact definitions.)
We shall show in the following that {\em the Miura
transformation defines a holomorphic foliation on a certain
complex manifold $X$.\/} More precisely, let
$\pi:X\rightarrow C$ be the relative complete flag $Fl(E),$
for $(E,\nabla )$ the flat holomorphic vector bundle on $C$
introduced before.  One can define (relative) Schubert Cells
in $X$ in the same fashion as in the standard case.
Moreover, there is a natural holomorphic foliation on $X$
induced by the horizontal distribution corresponding to
$\nabla.$

Then {\em the Miura transformation is the solution to
  the differential system defining this foliation in the
  affine coordinates on the ``Big Cell'' $X_0\subset X.$
  ($X_0$ is an affine bundle over $C.$)\/} Moreover, the
{\em operator \eqref{DS-form} arises as a partial connection
  $\delta$ associated to the same foliation on $X.$ By
  construction, the pull-back connection on $\pi^*E$
  obtained from $\nabla$ is adapted to the partial
  connection $\delta.$\/}

\subsection{}
The organization of this paper will be as follows. In the next
section we set up the notation and describe the various
construction we need in the following. Next, we discuss the
holomorphic foliation and describe its relation to the Miura
transformation. Finally, we briefly discuss some examples
and draw our conclusions.

\section{Some constructions}
\subsection{Notation and setup}
The compact Riemann surface $C$ of genus $g>1$ will be kept
fixed throughout. Local coordinates on $C$ will customarily
be denoted by $z,z_\alpha ,z_\beta ,\dots $ with domains
$U,U_\alpha ,U_\beta ,\dots$ All bundles and morphisms will
be considered in the holomorphic category.

\subsubsection{}
For the sake of simplicity we
shall restrict ourselves in the following to the special
linear group.
We shall use the canonical root decomposition
$\frak{g}\equiv\liesln =
\frak{n}^+\oplus\frak{h}\oplus\frak{n}^-$ such that the
nilpotent subalgebras $\frak{n}^\pm$ correspond to strictly
upper (lower) tringular matrices and the Cartan subalgebra
$\frak{h}$ to diagonal ones.
$\frak{b}^\pm=\frak{h}\oplus\frak{n}^\pm$ are the upper
(lower) Borel subalgebras and $N^\pm,\, B^\pm$ are the
analytic subgroups corresponding to the subalgebras so far
introduced.

\subsubsection{}\label{notation1}
We consider the pair $({E},\nabla )$ where ${E}\rightarrow
C$ is the jet-bundle mentioned before. It is determined by
the following condition~\cite{AlFa}. Let us temporarily
append the subscript ``$n$'', so that $E_{(n)}\equiv E;$
then $E_n$ is characterized as the extension ($E_1\equiv
{\script O}_C$)
\begin{displaymath}
  0\longrightarrow K_C^{\frac{1}{2}}\otimes
  E_{(n-1)}\longrightarrow E_{(n)}\longrightarrow
  K_C^{-\frac{n-1}{2}}\longrightarrow 0
\end{displaymath}
It follows that it is equipped with a complete filtration by
subbundles
\begin{equation}
  \{ 0\}\equiv {F}^0_{(n)}\subset
  {F}^1_{(n)}\subset\dots\subset {F}^{n-1}_{(n)}\subset
  {F}^n_{(n)}\equiv E_{(n)}
\label{filtration}
\end{equation}
with $F^r_{(n)}\cong K_C^{\frac{n-r}{2}}\otimes E_{(r)}$ and
${F}^{r+1}/{F}^r\cong {{K}_C}^{\frac{n-1}{2}-r}@, ,\,
r=0,\dots ,n-1,$ so that $\mathrm{Gr}F^\bullet\cong
\widetilde{E}$.  Explicit expressions for the cocycle
$\{\varphi_{\alpha\beta}\}$ representing ${E}$ can be found
in~\cite{AlFa,Zucc}; what is relevant for us is that they
can be put into {\em lower triangular form.\/} Thus there
are isomorphisms ${E}|_U\rightarrow {\script O}_C^n|_U$
together with bases $e_1,\dots ,e_n$ of sections of ${E}$ on
$U\subset C$ such that
${F}^r|_U$ is identified with the span $\langle
e_{n+1-r},\dots ,e_n\rangle$. $\nabla$ is a flat analytic
connection whose local expression is given by \eqref{W-form}
above. According to that, we have $\nabla {F}^r
\subset {F}^{r+1}\otimes{K}_C$.

\subsubsection{}
$\pi :Fl(E)\rightarrow C$ is the relative flag whose fiber
over a point $p\in C$ is the manifold $Fl(E_p)$ of complete
flags in $E_p$. We shall simply set $X\equiv Fl({E}).$ On
$Fl(E)$ we have the usual universal (tautological) flag $\{
0\}\subset S^{1}\subset\dots\subset S^{n-1}\subset \pi^*{E}$
and a corresponding universal sequence of quotients
$\pi^*E\twoheadrightarrow
Q^{n-1}\twoheadrightarrow\dots\twoheadrightarrow Q^{1},$
where $Q^{n-i}\equiv \pi^*E/S^{i}.$ Recall that if
$W^\bullet\in X$ then $S^i_{W^\bullet}=W^i.$

\subsubsection{}
We set $\rho :{P}\rightarrow C$ to be the $\sln$ principal
bundle corresponding to ${E}$. (Since
$\det{E}\cong\script{O}_C$, the structure group is reduced
to $\sln$.) Furthermore, restricting the action of $\sln$ on
${P}$ to $B^+$ yields
\begin{displaymath}
\begin{CD}
{P}@>{\tilde\pi}>>X\\
@V{\rho}VV               @VV{\pi}V\\
C       @=              C
\end{CD}
\end{displaymath}
as a relative version of the natural projection $\sln
\overset{B^+}{\longrightarrow} Fl(\CC ^n)$.
\begin{remark}
  Notice that
  $B^+\hookrightarrow{P}\overset{{\tilde\pi}}{\rightarrow}X$
  {\em does not\/} correspond to the reduction given by the
  filtration \eqref{filtration}. In fact the latter reduces
  the structure group to the {\em lower\/} Borel $B^-$.
\end{remark}

\subsection{}
The flat connection $\nabla$ determines a horizontal
holomorphic integrable distribution ${\script
  H}\hookrightarrow T_{ P}.$ The image ${\script
  F}\equiv\tilde\pi_*{\script H}$ of the horizontal
distribution ${\script H}$ is clearly holomorphic and
integrable (both by dimensional reasons and because it is
the direct image of an integrable distribution), so it
defines a holomorphic foliation in $X.$ Moreover, the
integrable distribution ${\script H}$ in $P$ satisfies with
respect to the fibration
$P\overset{{\tilde\pi}}{\rightarrow} X$ all the properties
of the horizontal distribution of a connection, except that
it does not project onto the tangent bundle of the base but
only on a foliation thereof, a situation referred to as {\em
  partial connection.\/}~\cite{KaTo}

The datum of a partial connection determines an operator on
$X$
\begin{displaymath}
\delta:\pi^*{E}\longrightarrow\pi^*{E}\otimes{\script F}^*
\end{displaymath}
satisfying the Leibnitz rule in the form
\begin{displaymath}
\delta (fs) = \tau (d f)\otimes s + f\delta (s)
\end{displaymath}
where $\tau$ is the projection $\tau:\Omega^1_X\rightarrow
{\script F}^*$ resulting from the injection ${\script
F}\hookrightarrow T_X.$

The connection $\nabla$ pulls back to a flat and holomorphic
connection on $\pi^*E$ (and $\pi^*P$), denoted by the same
symbol. This pull-back connection is {\em adapted\/} to the
partial connection if $\delta = (1\otimes\tau )\circ\nabla
.$ (see ref.~\cite{KaTo} and below.)

\begin{remark}
  It is interesting to see what the leaves of ${\script F}$
  are. Since $E$ is flat, $X\cong
  \widetilde{C}\times_{\pi_1(C)} Fl(\CC^n),$ where
  $\widetilde{C}$ is the universal cover of $C$ and
  $\pi_1(C)$ acts on the flag $ Fl(\CC^n)$ through the
  holonomy representation. Then the leaves are just the
  images of $\widetilde{C}\times \{{\mathrm p}{\mathrm t}\}$
  under the projection $\widetilde{C}\times
  Fl(\CC^n)\rightarrow \widetilde{C}\times_{\pi_1(C)}
  Fl(\CC^n).$
\end{remark}

\subsection{}
We set out to characterize geometrically the calculations
presented in \ref{local-1}, \ref{local-2}. Following ref.
\cite{Ful92}, we use the filtration \eqref{filtration} to
define a submanifold $X_0\subset X.$ As a set, it is
determined by the condition that for any $p\in C$ it
consists of those flags $W^\bullet :\{0\}\subset
W^1\subset\dots\subset W^{n-1}\in X_p=\pi^{-1}(p)$ such that
for any $i$
\begin{equation}\label{big-cell}
  \dim (W^{j}\cap F^{n+i-j}_p)=i\, .
\end{equation}
It is easy to see that this condition is in fact redundant
and the very same locus is characterized by restricting
\eqref{big-cell} to $i=0$ only. More invariantly, one is led
to the equivalent condition that $X_0$ is the locus where
the map
\begin{displaymath}
  \pi^*F^k\hookrightarrow\pi^*E\twoheadrightarrow Q^k
\end{displaymath}
has maximal rank \cite{Ful92}.
We can condense its geometrical meaning in the following
\begin{lemma}
  $X_0$ is a relative big cell in $X$, i.e. for any $p\in
  C,$ $X_0\cap X_p$ is the big cell in $X_p\cong Fl(\CC^n)$,
  namely the orbit of the lower unipotent subgroup of $\sln$
  through the standard flag. Furthermore,
  $X_0\overset{\pi|_{X_0}}{\longrightarrow}C$ is an affine
  bundle over $C.$
\end{lemma}
\subsubsection{}
We sketch the proof in order to set some convention needed
in the following.  Full details can be found in
\cite{Ful92}. It follows from \ref{notation1} that we can
use coordinates in each fiber such that:
\begin{enumerate}
\item The canonical projection $\sln\rightarrow Fl(\CC^n)$
is realized by sending $g\in\sln$ to the flag $\langle
g\cdot e_1\rangle\subset\langle g\cdot e_1,g\cdot
e_2\rangle\subset\dots\subset\CC^n$ where $\langle
e_1\rangle\subset\langle e_1,e_2\rangle\subset\dots$ is the
standard flag, whose stabilizer is the upper Borel $B^+.$
\item The filtration corresponds to the flag
\begin{displaymath}
V^\bullet:V^1\subset\dots\subset V^{n-1}\subset V^n\equiv
\CC^n ,
\end{displaymath}
where $V^i=\langle e_{n+1-i},\dots ,e_n\rangle,$ and it comes
from the longest permutation $e_i\rightarrow e_{n+1-i}$ (up
to a sign).
\end{enumerate}
In this way we are reduced to the non relative situation.
Geometrically, a flag $W^\bullet$ in $Fl(\CC^n)$ satisfies
condition \eqref{big-cell} if the fixed flag $V^\bullet$
sweeps a complete flag in each subspace $W^{j}.$ On the
other hand, it is just a computation to verify that in this
case $W^\bullet$ is obtained from an element in $\sln$ {\em
  all of whose principal minors are non singular,\/} that is
an element in the big cell of $\sln.$ This element can
always be factorized into a product of a lower unipotent
matrix times an upper Borel one. It follows that the flag
$W^\bullet$ is represented by that lower unipotent matrix,
whose entries together with $z$ are thus coordinates on
$X_0.$ Furthermore, let $p\in U_\alpha\cap U_\beta\subset C$
and let $W^\bullet\in X_0$ be a flag over $p.$ It follows
that $W^\bullet$ can be represented by a lower triangular
matrix (frame) $\bold{N}_{-,\alpha}$ with respect to any
trivialization and it is immediate that
$\bold{N}_{-,\alpha}\varphi_{\alpha\beta}^{0}(p)=
\varphi_{\alpha\beta}(p)\bold{N}_{-,\beta}$ in the overlap,
where $\varphi_{\alpha\beta}^{0}$ is the diagonal part of
$\varphi_{\alpha\beta}.$ This relation shows that $X_0$ is
an affine bundle over $C.$

\subsubsection{}\label{frames}
It follows from \eqref{big-cell} and the sketch of the proof
above that $E|_{X_0}$ becomes
(holomorphically) decomposable, hence we get at once the
\begin{corollary}
\begin{displaymath}
\left. \pi^* E\right|_{X_0}\cong
  \left. \pi^*\widetilde{E}\right|_{X_0}\qed
\end{displaymath}
\end{corollary}
Thus there are {\em frames\/} $\bold{M}$ in $P$ at $p\in C$
---~those projecting down to $X_0$~--- that can be
factorized into $\bold{M} =\bold{N}_-\bold{B}_+,$ formally
as in \eqref{gauss}, and the entries of $\bold{N}_-$ are
coordinates on $X_0\cap X_p.$ Moreover, $\bold{N}_-$ can be
considered as a new frame in $\pi^*E$ which is ``adapted''
to the tautological flag $S^\bullet\hookrightarrow\pi^*E.$

\section{Miura transformation and the horizontal foliation}
\subsection{}
In order to characterize the distribution ${\script F}$ on
$X,$ recall that ${\script F}\equiv \pi_*{\script H},$ where
${\script H}$ is the horizontal distribution in $P$
determined by the connection $\nabla .$ Its ``equations''
are simply given by $\omega =0,$ where $\omega$ is the
corresponding connection form. In particular we get
$\omega^0=\omega^+=\omega^-=0,$ where
$\omega^0,\,\omega^\pm$ are the components of $\omega$ with
respect to the root decomposition of the Lie algebra
$\liesln .$
\begin{lemma}
  The form $\omega^-$ is horizontal with respect to the
  fibration $P\overset{{\tilde\pi}}{\rightarrow}X,$ hence it
  descends on $X.$ The equation for ${\script F}$ in $X$ is
  $\omega^-=0.$
\end{lemma}
\begin{pf}
  The relative flag $X$ is gotten as a quotient of $P$ by
  the action of $B^+,$ so that we get $\omega^-(\xi )=0$ for
  any $\xi\in\frak{b}^+.$ (We identify elements of the Lie
  algebra with the corresponding fundamental vector fields.)
  Since they generate the vertical bundle of
  $P\overset{{\tilde\pi}}{\rightarrow}X,$ the
  $\frak{n}^-$--valued 1--form $\omega^-$ is {\em
    horizontal\/} and we can keep the same name for the
  corresponding object on $X.$ It follows that an element of
  $T_X$ is in ${\script F}$ iff it annihilates $\omega^-.$
\end{pf}

\subsection{}\label{coordinates}
We now introduce explicit coordinates and restrict our
attention to $X_0\subset X.$ If $p\in U\subset C,$ with $z$
a local parameter at $p,$ we trivialize $P|_U\cong
U\times\sln$ with coordinates $(z,g).$ Moreover, if
$(z,g)\in P|_U\cap {\tilde\pi}^{-1}(X_0)$ then $g=n_-b_+,$
with $b_+\in B^+ ,\, n_-\in N_-.$ Thus ${\tilde\pi}:P|_U\cap
{\tilde\pi}^{-1}(X_0)\rightarrow X_0|_U$ is locally given by
$(z,g)\mapsto (z,n_-).$ Furthermore, let
$\sigma_-:\AAA^{n(n-1)/2}\rightarrow N_-$ be the map sending
an $n(n-1)/2$-tuple of complex numbers into a lower
unipotent matrix. Since $(z,\underline{u})\in
U\times\AAA^{n(n-1)/2}$ are local coordinates on $X_0|_U,$
$(z,\underline{u})\mapsto (z,\sigma_-(\underline{u}))$ is a
local section to $P\overset{{\tilde\pi}}{\rightarrow}X.$ It
follows that
\begin{align}
  U\times\AAA^{n(n-1)/2}\times B^+ &\longrightarrow P|_U\cap
  {\tilde\pi}^{-1}(X_0)\notag\\ (z,\underline{u},b) &\mapsto
  (z,\sigma_-(\underline{u})\, b)\notag
\end{align}
is a new choice of local coordinates on $P.$ Putting
\begin{equation}\label{splitting}
  g = \sigma_-(\underline{u})\, b
\end{equation}
is a coordinate change in the fibers of $P$ over points in
$X_0|_U\subset X.$

\subsection{}
Let $A$ be the local connection matrix for $\nabla$ on $U$
as given in \eqref{W-form}, so that
\begin{equation}
  \omega|_{P|_U}=g^{-1}dg + \ad_{g^{-1}}(A)
\end{equation}\label{omega-local}
is the local expression for the connection form, where
$g^{-1}dg$ denotes the Maurer-Cartan form on $\sln.$ Using
\eqref{splitting} we compute
\begin{equation}\label{omega-local-new}
  \omega|_{P|_U}=b^{-1}db +
  \ad_{b^{-1}}(\sigma_-^{-1}d\sigma_- +
  \ad_{\sigma_-^{-1}}(A))
\end{equation}
and since $b\in B^+$ and $\sigma_-(\underline{u})\in N_-$ we
have
\begin{equation}\label{omega-meno-local} \omega^-|_{P|_U} =
  \sigma_-^{-1}d\sigma_- + {(\ad_{\sigma_-^{-1}}(A))}^-
\end{equation}
where $(\cdot)^-$ denotes the projection from $\liesln$ to
$\frak{n}^-.$ Thus the differential ideal determining
${\script F}$ on $X_0$ is generated by the entries of the
$\frak{n}^-$-valued 1-form $\omega^-|_{P|_U}$
\eqref{omega-meno-local} above.

The quantity in parentheses in \eqref{omega-local-new} is
interesting, insofar it is formally analogous to
\eqref{gauge}. Indeed, let us explicitly write
$\nabla\bold{e} = \bold{e}\otimes A$ for $\nabla$ on $E$ and
$\bold{e}=(e_1,\dots ,e_n)$ the frame introduced in
\eqref{notation1}. Let us also keep the same notation for
the corresponding object pulled back to $X|_U$ under $\pi .$
{}From \ref{frames} and \ref{coordinates} we have that
$\bold{e}\cdot\sigma_-(\underline{u})$ is a new frame for
$\pi^*E\cong\pi^*\widetilde{E}$ on $X_0|_U,$ and the new connection
matrix on $X_0|_U$ reads
\begin{align}
  {\tilde A} &= \sigma_-^{-1}d\sigma_- +
  \ad_{\sigma_-^{-1}}(A)\notag\\ &=
\begin{pmatrix}
  u_1&1&0&\dots&0\\ 0&u_2-u_1&1&\dots&0\\
  \dots&&\dots&&\dots\\ \hdotsfor{3}&u_{n-1}-u_{n-2}&1\\
  0&\dots&\dots&0&-u_{n-1}
\end{pmatrix}dz + \omega^-|_{P|_U}
\end{align}
Now, $\tau :\Omega^1_X\rightarrow{\script F}^*$ kills the
annihilator of ${\script F}$ and hence $\omega^-,$ therefore
$\tau ({\tilde A})$ has the same form as the ``connection
matrix'' in \eqref{DS-form}. It follows that the diagonal
gauge operator is best interpreted as the partial connection
operator $\delta$ on $X_0|_U.$ Note that by construction the
connection $\nabla$ on $\pi^*E$ is adapted to $\delta ,$
which is what we wanted to show.

\section{Examples}
\subsection{Rank 2} The bundle $E$ is the well known
extension~\cite{Gun66}
\begin{displaymath}
0\longrightarrow K_C^{1/2}\longrightarrow E\longrightarrow
K_C^{-1/2}\longrightarrow 0
\end{displaymath}
with transition functions
\begin{displaymath}
  \varphi_{\alpha\beta}=\begin{pmatrix}
  k_{\alpha\beta}^{-1/2}&0\\
  \frac{d}{dz_\beta}k_{\alpha\beta}^{1/2}
  &k_{\alpha\beta}^{1/2} \end{pmatrix}
\end{displaymath}
where $k_{\alpha\beta}=(dz_\alpha /dz_\beta)^{-1}$ are the
transition functions for $K_C.$ The connection $\nabla$ has
the form (see \eqref{W-form})
\begin{displaymath}
  \nabla = \partial +
\begin{pmatrix}
0&1\\T(z)&0
\end{pmatrix}dz\, .
\end{displaymath}
In this case $X=\PP(E)$ and $X_0=X\setminus C_\infty,$ where
$C_\infty$ is the divisor at infinity, that is the image of
$C$ under the canonical section of the ruled surface $X$
given by the sub-linebundle $K_C^{1/2}.$ Points in $X_0$ are
of the form $\left[\begin{smallmatrix} 1\\u
 \end{smallmatrix}\right]$
and the section
$\sigma_-:X_0|_U\rightarrow P|_{X_0|U}$ is given by
$\sigma_-(u)=\left(
\begin{smallmatrix}
1&0\\u&1
\end{smallmatrix}\right).$ It
follows then that the 1-form defining the foliation on
$X_0\subset X$ has the expression
\begin{equation}\label{miura-fol-2}
\omega^-|_{P|_U} = dw + (T(z) - w^2)dz
\end{equation}
which should be compared with expression \eqref{miura-2}.

\subsection{Rank 3} In this case $E=J^2(K_C^{-1}).$ An
explicit cocycle for it is~\cite{AlFa}
\begin{displaymath}
\varphi_{\alpha\beta}=
\begin{pmatrix}
k_{\alpha\beta}^{-1}&0&0\\
2\sigma_{\alpha\beta}&1&0\\
\partial^2k_{\alpha\beta}&2
k_{\alpha\beta}\sigma_{\alpha\beta}&k_{\alpha\beta}
\end{pmatrix}
\end{displaymath}
where $k_{\alpha\beta}$ has the same meaning as before and
$\sigma_{\alpha\beta}=\partial\log k_{\alpha\beta}^{1/2}.$

It follows from the general theory that $X_0$ consists of
flags $W^1\subset W^2$ such that $S^2_{W^\bullet}\cap
(\pi^*F^1)_{W^\bullet}=\{ 0\}$ and $S^1_{W^\bullet}\cap
\big(S^2_{W^\bullet}\cap (\pi^*F^2)_{W^\bullet}\big)=\{
0\}.$ Local coordinates are $(z,u_1,u_2,u_3)$ and
$\sigma_-=\left(
\begin{smallmatrix}
1&0&0\\u_1&1&0\\u_3&u_2&1
\end{smallmatrix}\right)$ is the section of
$P\overset{\tilde\pi}{\rightarrow}X$ on $X_0|_U.$ The
foliation is given by the following differential forms
\begin{align*}
\omega^-_1 &= du_1 + (u_3 - u_1^2)\, dz\\
\omega^-_2 &= du_2 + (T(z) - u_2^2 + u_1 u_2 - u_3)\, dz\\
\omega^-_3 &= du_3 - u_2 du_1 +
(w(z) + u_1 T(z) - u_2 u_3 + u_1^2 - u_1 u_3)\, dz
\end{align*}
and it is easily verified that the solution of this
differential system is precisely \eqref{miura-3}.

\section{Conclusions}
We have pointed out how the global geometry underlying the
Miura transformation is to be found in the foliation
${\script F}$ in the flag bundle $X=Fl(E).$ As mentioned in
the introduction, one of the problems concerning the
definition of the Miura transformation on a Riemann surface
was the monodromy behavior of its solutions.  According to
the point of view outlined here, this monodromy is nothing
but the way the leaves are stacked in the foliation, thus it
is an essential part of the overall structure. Since
${\script F}$ comes from a flat structure, it is moreover
just another way to look at the holonomy of the flat
connection. The next step would be to set up a proper
generalization of the free field quantization scheme in
terms of this new interpretation of the Miura
transformation. We hope to return to the subject elsewhere.

\subsection*{Acknowledgements} I would like to thank
J. L. Dupont, G. Falqui and F. Kamber for discussions, as
well as L. Bonora and L. Takhtajan for a careful reading of
the manuscript.



\begin{thebibliography}{10}
\bibitem{AlBo}
E. Aldrovandi, L. Bonora, {\em Liouville and Toda Field
Theories on Riemann Surfaces.\/} \jgp{14}{1994}{65--109},
and \hepth{9303064}
\bibitem{ABBP}
E. Aldrovandi, L. Bonora, V. Bonservizi and R. Paunov, {\em
Free Field representation of Toda Field Theories.\/}
\ijmp{9}{1994}{57--86}, and \hepth{9211112}
\bibitem{AlFa}
E. Aldrovandi, G. Falqui, {\em Geometry of
Higgs and Toda Fields on Riemann Surfaces\/.} \hepth{9312093},
to appear in {\em Journal of Geometry and Physics.\/}
\bibitem{Atiyah}
M. F. Atiyah, {\em Complex analytic connections in fibre
bundles.\/} \tams{85}{1957}{185--207}.
\bibitem{BiGe}
A. Bilal, J.--L. Gervais,
{\em Extended $c=\infty$ Conformal Systems from classical
  Toda field theories,\/} \np{314}{1989}{646--686}; {\em
  Systematic construction of Conformal Field Theories with
  higher--spin Virasoro algebra,\/} \np{318}{1989}{579--630}.
\bibitem{BoX}
L. Bonora, C. S. Xiong, {\em Covariant $sl_2$
decomposition of the $sl_n$ Drinfel'd--Sokolov Equations ad
the $W_n$-algebras,\/} \ijmp{7}{1992}{1507--1525}.
\bibitem{ItZu}
P. DiFrancesco, C. Itzykson, J.B. Zuber, {\em Classical
  $W$-algebras,\/} \cmp{140}{1991}{543--568}.
\bibitem{FaLy} V. A. Fateev, S. L. Lykyanov, {\em The Models
of two-dimensional Conformal Quantum Filed Theory with $Z_n$
symmetry,\/} \ijmp{3}{1988}{507--520}.
\bibitem{Fre95}
E. Frenkel, {\em Affine Algebras, Langlands Duality and
  Bethe Ansatz,\/} \qalg{9506003}.
\bibitem{Ful92}
W. Fulton, {\em Flags, Schubert polynomials, Degeneracy
loci, and determinantal formulas.\/} \dmj{65}{1992}{381--420}.
\bibitem{GeNe}
J.--L., Gervais, A. Neveu,
{\it The dual string spectrum in
Polyakov's quantization. I.\/} \np{199}{1982}{59};
{\it The dual string spectrum in
Polyakov's quantization. II.\/} \np{209}{1982}{125}.
\bibitem{Gun66}
R. C. Gunning,
{\em Lectures on Riemann Surfaces.\/} Princet. Univ. Press.,
Princeton 1966.
\bibitem{Gun67}
R. C. Gunning,
{\em Lectures on Vector bundles over Riemann Surfaces.\/}
Princeton. Univ. Press., Princeton 1967.
\bibitem{KaTo}
F. Kamber, Ph. Tondeur,
{\em Foliated Bundles and Characteristic Classes.\/}
Lecture Notes in Math. {\bf 493}, Springer, 1975.
\bibitem{suresh}
S. Govindarajan, {\em Higher Dimensional Uniformization and
$W$-geometry,\/} \hepth{9412078}.
\bibitem{Zucc}
R. Zucchini, {\em Light cone $W_n$ geometry and its
symmetries and projective field theory.\/}
Class. Quant. Grav. {\bf 10} (1993), 253--278, and
\hepth{9205102}
\bibitem{ZuDS}
R. Zucchini, {\em A Krichever--Novikov
formulation of $W$--algebras on Riemann Surfaces,\/}
\pl{323}{1994}{322--329}, and \hepth{9310061}
\end{thebibliography}
\end{document}